\documentclass[12pt,twoside]{article}
\usepackage{amsmath}
\usepackage{amsfonts}
\usepackage{amsthm}
\usepackage{amssymb}
\usepackage{graphicx}
\usepackage{multicol}

\begin{document}
\title{{\normalsize{\bf Kervaire conjecture on weight of group via fundamental group of 
ribbon sphere-link}}}
\author{{\footnotesize Akio KAWAUCHI}\\
{\footnotesize{\it Osaka Central Advanced Mathematical Institute, Osaka Metropolitan University }}\\
{\footnotesize{\it Sugimoto, Sumiyoshi-ku, Osaka 558-8585, Japan}}\\
{\footnotesize{\it kawauchi@omu.ac.jp}}}
\date\, 
\maketitle
\vspace{0.25in}
\baselineskip=15pt
\newtheorem{Theorem}{Theorem}[section]
\newtheorem{Conjecture}[Theorem]{Conjecture}
\newtheorem{Lemma}[Theorem]{Lemma}
\newtheorem{Sublemma}[Theorem]{Sublemma}
\newtheorem{Proposition}[Theorem]{Proposition}
\newtheorem{Corollary}[Theorem]{Corollary}
\newtheorem{Claim}[Theorem]{Claim}
\newtheorem{Definition}[Theorem]{Definition}
\newtheorem{Example}[Theorem]{Example}

\begin{abstract} 
Kervaire conjecture that the weight of the free product of every non-trivial group and the infinite cyclic group is not one is affirmatively  confirmed by confirming affirmatively 
Conjecture~${\mathbf Z}$ on the knot exterior   introduced by  Gonz{\'a}lez Acu{\~n}a and 
Ram{\'i}rez   as a   conjecture equivalent to  Kervaire conjecture.  

\phantom{x}



\end{abstract}

\maketitle

\section*{1. Introduction} 

A {\it weight system} of a group $G$ is a system of elements $w_i,\,(i=1,2,\dots,n)$ of $G$ 
such that the normal closure $N(w_1,w_2,\dots, w_n)$ of $w_i,\,(i=1,2,\dots,n)$ 
in $G$ (=: the smallest normal subgroup generated by $w_i,\,(i=1,2,\dots,n)$ in $G$) is equal to $G$. 
The {\it weight} of a group 
$G$ is the least cardinal number $w(G)$ of a weight system of $G$. By convention, 
 $w(G)=0$ if and only if $G$ is the trivial group. 
The {\it rank} of $G$ is the least cardinal number $r(G)$ of generators of $G$. The  difference $r(G)- w(G)$ is non-negative and  in general taken sufficiently large. For example, 
let $G=\pi_1(S^3\setminus k, x)$ be the fundamental group of a polygonal knot $k$ in $S^3$.
Then $G \cong {\mathbf Z}$ and $r(G)=1$ for the trivial knot $k$, 
 $r(G)=2$  for the trefoil knot $k=3_1$, and 
$r(G)=n$ for the $n-1(\geq 2)$-fold connected sum $k=\#_{n-1} 3_1$ of the trefoil knot $3_1$.
On the other hand, $w(G)=1$ for every knot $k$, 
because $G/N(m(k))=\{1\}$ for a meridian element $m(k)$ of $k$. 
Let $G*{\mathbf Z}$ denote the free product  of  a group $G$ and the infinite cyclic group 
${\mathbf Z}$. 
{\it Kervaire's conjecture} on the weight of a group is the following conjecture (see Kervaire \cite{Ker}, Magnus-Karrass-Solitar \cite[p. 403]{MKS}): 

\phantom{x}

\noindent{\bf Kervaire Conjecture.} 
$w(G*{\mathbf Z})>1$ for   every non-trivial group $G$. 

\phantom{x}

Some partial affirmative confirmations of this conjecture are known. For example, the following result of  Klyachko \cite{Klyachko} is used in this paper: 

\phantom{x}

\noindent{\bf Theorem (Klyachko).}  $w(G*{\mathbf Z})>1$ for every 
non-trivial torsion-free group $G$.

\phantom{x}

A {\it knot exterior} is a compact 3-manifold $E=\mbox{cl}(S^3\setminus N(k))$ for a tubular neighborhood $N(k)$ of  a polygonal knot $k$ in the 3-sphere $S^3$. 
Let  $F$ be a compact connected orientable  non-separating proper surface of  $E$ where 
the boundary $\partial F$ of $F$ may be disconnected. 
Let $E(F)=\mbox{cl}(E\setminus F\times I)$ be the compact 
piecewise-linear 3-manifold for a normal line bundle $F\times I$ of $F$ in $E(F)$ 
where $ I=[-1,1]$. 
Let $E(F)^+$ be the 3-complex obtained from $E(F)$ by adding the cone 
$\mbox{Cone}(v,F\times \partial I)$  over the base $F\times \partial I$ with a  vertex $v$ 
disjoint from $E$, where $\partial I=\{ 1,-1\}$. 
The 3-complex $E(F)^+$ is also considered to be  obtained from $E$ by shrinking the  normal line bundle $F\times I$ into the vertex $v$.  
The result of Conjecture~${\mathbf Z}$ due to Gonz{\'a}lez Acu{\~n}a and Ram{\'i}rez in \cite{GoR} is stated as follows:

\phantom{x}

\noindent{\bf Theorem (Gonz{\'a}lez Acu{\~n}a-Ram{\'i}rez).} 
Kervaire's conjecture is equivalent to the following  conjecture: 

\medskip 

\noindent{\bf Conjecture~${\mathbf Z}$.}  The fundamental group 
$\pi_1(E(F)^+,v)$  is isomorphic to ${\mathbf Z}$ for every knot exterior $E$ and every 
compact connected orientable  non-separating proper surface $F$ in $E$.  

\phantom{x}

See \cite{EM,GoR,RG} for some knot theoretical investigations of this surface $F$ and some partial confirmations. 
In this paper, Kervaire conjecture is  confirmed affirmatively by confirming 
Conjecture~${\mathbf Z}$ affirmatively. 

\phantom{x}

\noindent{\bf Theorem~1.} 
Conjecture~${\mathbf Z}$ is true. 

\phantom{x}

Gonz{\'a}lez Acu{\~n}a-Ram{\'i}rez theorem and Theorem~1 imply: 

\phantom{x}

\noindent{\bf Corollary~2.} 
Kervaire conjectureis true. 

\phantom{x}

An outline of the proof of  Theorem~1 is explained as follows.

\phantom{x}

\noindent{\bf Outline of the proof of  Theorem~1.} Let 
\[E(F)^{++}=E(F) \cup \mbox{Cone}(v_+,F\times 1)\cup \mbox{Cone}(v_-,F\times (-1))\]
be a 3-complex for distinct vertexes $v_+$ and $v_-$ disjoint from $E$. 
Then the 3-complex $E(F)^+$ is homotopy equivalent to the bouquet 
$E(F)^{++}\vee S^1$.
Hence  the fundamental group $\pi_1(E(F)^+,v)$ is isomorphic to the free 
product $\pi_1(E(F)^{++},v)*{\mathbf Z}$.
Thus, $\pi_1(E(F)^+,v) \cong  {\mathbf Z}$ if and only if  $\pi_1(E(F)^{++},v) =\{1\}$ 
and Conjecture~${\mathbf Z}$ is equivalent to the claim that  $\pi_1(E(F)^{++},v) =\{1\}$. 
The following observation is used.

\phantom{x}

\noindent{\bf Lemma~3.}
$w(\pi_1(E(F)^+,v))= w(\pi_1(E(F)^{++},v)*{\mathbf Z})=1$. 

\phantom{x}

\noindent{\bf Proof of  Lemma~3.} 
Because the fundamental group $\pi_1(E(F)^+,v)$ is a non-trivial 
quotient group of   $\pi_1(E,v)$ and $w(\pi_1(E,v))=1$, 
the desired result is obtained.
This completes the proof of Lemma~3. $\square$.

\phantom{x}

The following lemma is proved in Section~2.

\phantom{x}

\noindent{\bf Lemma~4.}
The fundamental group $\pi_1(E(F)^+,v)$ is a torsion-free group. 

\phantom{x}

By assuming Lemma~4, the proof of Theorem~1 is completed as follows:

\phantom{x}

\noindent{\bf Proof of Theorem~1.} 
Klyachko Theorem says that 
if $G$ is a torsion-free group and $w(G*{\mathbf Z})=1$, then $G=\{1\}$. 
Hence by this theorem and  Lemmas~3, 4, $\pi_1(E(F)^{++},v)\cong\{1\}$ and 
$\pi_1(E(F)^+,v)\cong {\mathbf Z}$. 
This completes the proof of Theorem~1. $\square$.

\phantom{x}

In the first draft of this research, the author tried to show that every finitely presented group 
$G$ with  $w(G*{\mathbf Z})=1$  is torsion-free. 
This trial succeeds for a group $G$ 
of deficiency $0$, but failed for a group $G$ of negative deficiency. 
The main point of this failure is the attempt to construct a finitely presented group of deficiency $0$ from the group of negative deficiency, which  forced the author to show that $G$ is 
a torsion-free group while the deficiency remains negative.
Fortunately, the fundamental group $\pi_1(E(F)^+,v)$ of the 3-complex $E(F)^+$  was an excellent object to this consideration, so it could be done. 

\phantom{x}

\section*{2. Proof of Lemma~4 } 

The proof of  Lemma~4 is done  as follows by using the consept of collapse in \cite{Hud}. 

\phantom{x}

\noindent{\bf Proof of Lemma~4.} 
Collapse $F$ into a triangulated graph $\gamma$ by using that  $F$ is a bounded surface. 
Enlarge the fiber $I$  of a normal line bundle $F\times I$  of $F$ in $E$ into a fiber $J$ 
of a normal line bundle $F\times J$ of $F$ in $E$
so that  $I\subset  J\setminus \partial J$.   
Let $J^c=\mbox{cl}(J \setminus I)$.
Let $E(F)^-=\mbox{cl}(E\setminus F\times J)$. 
Collapse $F\times J^c$ into $\gamma\times J^c$. 
Triangulate $\gamma\times J^c$ without introducing new vertexes. 
The 3-complex  $E(F)^+$ is collapsed into  a finite 3-complex 
\[E(F)^-\cup \gamma\times J^c \cup \mbox{Cone}(v,\gamma\times \partial I)\]
and thus collapsed into  a finite 2-complex 
\[ P=P^- \cup \gamma\times J^c \cup \mbox{Cone}(v,\gamma\times \partial I)\]
obtained by taking any 2-complex $P^-$ collapsed from $E(F)^-$. 
This 2-complex $P$ is a subcomplex of a 3-complex 
\[ Q= \mbox{Cone}(v,P^-\cup \gamma\times J^c).\]
Since every 2-complex of $\gamma\times J^c$ contains at most 
one 1-simplex of  $\gamma\times \partial I$,  
every 3-simplex of $\mbox{Cone}(v,\gamma\times J^c)$ contains at most 
one 2-simplex of  $\mbox{Cone}(v,\gamma\times \partial I)$. 
Collapse every 3-simplex of $\mbox{Cone}(v,\gamma\times J^c)$ from 
 a 2-face containing $v$ and not belonging to   
$\mbox{Cone}(v,\gamma\times \partial I)$. 
Then collapse every 3-simplex of $\mbox{Cone}(v,P^-)$ from any 2-face containing 
the vertex $v$. 
Thus, the 3-complex $Q$ is collapsed to a finite 2-complex $C$ containing 
the 2-complex $P$ as a subcomplex. Since $Q$ is  collapsed to the vertex $v$,  
$C$ is a finite contractible 2-complex. It is shown in \cite{K1}  that 
every connected subcomplex of a finite contractible 2-complex is 
aspherical (see also \cite{K2}).  
Since the fundamental group of a connected aspherical 
complex is  a torsion-free group,  the group $\pi_1(P,v)$ is a torsion-free group.  
Note that this torsion-freeness comes from the torsion-freeness of 
the fundamental group of a ribbon $S^2$-link in the 4-sphere  $S^4$ which 
corresponds bijectively to the fundamental group of a connected subcomplex of a finite contractible 2-complex, as it is discussed  in \cite{K1}.
Actually, the group $\pi_1(P,v)$ is a ribbon $S^2$-knot group,  for 
$H_1(P;{\mathbf Z})\cong {\mathbf Z}$. 
Since $\pi_1(E(F)^+,v)$ is isomorphic to $\pi_1(P,v)$, the group $\pi_1(E(F)^+,v)$ is 
a torsion-free group. 
This completes the proof of Lemma~4. 
$\square$

\phantom{x}

\section*{Acknowledgements}  
This work was partly supported by JSPS KAKENHI Grant Numbers JP19H01788, JP21H00978 and MEXT Promotion of Distinctive Joint Research Center Program JPMXP0723833165.


\begin{thebibliography}{99}

\bibitem{EM} M. Eudave-Mu{\~n}oz, On knots with icon surfaces, Osaka Journal of Mathematics, 50 (2013), 271–285.

\bibitem {GoR} F. Gonz{\'a}lez Acu{\~n}a and A. Ram{\'i}rez, A knot-theoretic equivalent of
the Kervaire conjecture, J. Knot Theory Ramifications, 15 (2006), 471–478. 

\bibitem{Hud} J. F. P. Hudson, Piecewise-linear topology, Benjamin (1969).

\bibitem{K1} A. Kawauchi, Ribbonness of Kervaire's sphere-link in homotopy 4-sphere and its consequences to 2-complexes. Available from: 
https://sites.google.com/view/ kawauchiwriting


\bibitem{K2} A. Kawauchi, Whitehead aspherical conjecture via ribbon sphere-link. Available from: 
https://sites.google.com/view/ kawauchiwriting

\bibitem{Ker} M. A. Kervaire, On higher dimensional knots, in: Differential and combinatorial topology, Princeton Math. Ser. 27 (1965), 105-119, Princeton Univ. Press.

\bibitem {Klyachko} Ant. A. Klyachko, A funny property of a sphere and equations over groups, 
Comm. Algebra 21 (1993), 2555–2575.

\bibitem{MKS} W. Magnus, A. Karrass and D. Solitar, Combinatorial group theory: 
Presentations of groups in terms of generators and relations, Interscience Publisheres (1966).


\bibitem{RG} J. Rodr{\'i}guez-Viorato and F. Gonz{\'a}lez Acu{\~n}a, On pretzel knots and
Conjecture~${\mathbf Z}$, Journal of Knot Theory and Its Ramifications, 25 (2016) 1650012.


\end{thebibliography}
\end{document}